\documentclass{elsart}
\usepackage{amssymb}
\begin{document}
\begin{frontmatter}
\title{Simplification Techniques for Maps in Simplicial Topology}
\author{Rocio Gonzalez--Diaz}
\author{Pedro Real}
\address{Dept. of Applied Math. I,
University of Seville, Spain}
\thanks{Partially
supported by the PAICYT research project
 FQM--296 from Junta de Andalucia.
Web address: {\em http://www.us.es/gtocoma}.}


\maketitle

\begin{abstract}
    This paper offers an algorithmic solution to the problem of
    obtaining ``economical'' formulae for some maps in
 Simplicial Topology, having, in principle, a high computational cost in their evaluation. In
 particular,  maps of
 this kind are
 used for defining
 cohomology operations at the cochain
level. As an example, we obtain explicit combinatorial descriptions of Steenrod $k$th
powers exclusively in terms of face operators.
\end{abstract}

\end{frontmatter}

\newcommand{\ra}{\rightarrow}
\newcommand{\la}{\leftarrow}
\newcommand{\lra}{\longrightarrow}
\newcommand{\lla}{\longleftarrow}
\newcommand{\x}{\bar{x}}
\newcommand{\Ra}{\Rightarrow}
\newcommand{\D}{\displaystyle}
\newcommand{\scst}{\scriptscriptstyle}
\newcommand{\scsz}{\scriptsize}
\newcommand{\pa}{\partial}
\newcommand{\m}{\bar{m}}
\newcommand{\ot}{\otimes}
\newcommand{\EML}{{\mbox{\it EML}}}
\newcommand{\AW}{{\mbox{\it AW}}}
\newcommand{\SHI}{{\mbox{\it SHI}}}
\newcommand{\ti}{\times}
\newcommand{\al}{\alpha}
\newcommand{\be}{\beta}
\newcommand{\Z}{{\bf Z}}
\newcommand{\Fp}{{\bf F}_p}
\newcommand{\Zp}{{\bf Z}_p}
\newcommand{\Zq}{{\bf Z}_q}
\newcommand{\Fr}{{\bf F}_2}
\newcommand{\Zr}{{\bf Z}_2}
\newcommand{\cF}{{\cal F}}
\newcommand{\cP}{{\cal P}}
\newcommand{\cC}{{\cal C}}
\newcommand{\cn}{C_*^{\scst N}}
\newcommand{\cm}{C_m^{\scst N}}
\newcommand{\ov}{\overline}

\section{Introduction}
   In this paper we deal with problems in the field of Combinatorial Topology. We work with
simplicial sets, which provide  combinatorial descriptions of
topological spaces. A simplicial set (see \cite{May67}) is a
graded set $K=\{K_q\}_{q\geq 0}$ whose $q$-dimensional ``building
blocks'' are $q$--simplices and whose ``mortar'' is face
($\partial_i:K_{q+1}\ra K_{q}$) and degeneracy ($s_i:K_q\ra
K_{q+1}$) operators. It is an elementary fact that any
composition of face and degeneracy operators of a simplicial set
$K$ can be expressed in the ``normalized'' form: $$s_{j_t}\cdots
s_{j_1}\partial_{i_1}\cdots \partial_{i_s},$$ where $j_t > \cdots
>j_1\geq 0$ and $i_s > \cdots > i_1 \geq 0$, due to certain
commutativity properties.
   Roughly speaking, we are interested here not only in ``normalizing'' some
   compositions of face and degeneracy operators, but also in determining which of them
   involve exclusively
   face operators. In particular, we simplify compositions that are used for defining
 important cohomology operations such as Steenrod squares \cite{Ste47}, Steenrod
 $k$th powers \cite{Ste52} or Adem secondary cohomology operations \cite{Ade52,Ade58}.
 In fact, from a simplicial viewpoint and taking into account that we deal with homological
   information  given in terms of explicit chain homotopy equivalences \cite{Rea00,Gon00}, the
   description of
invariants in Algebraic Topology can be reduced to the study of
compositions of certain specific maps given essentially in terms
of face and degeneracy operators. The fundamental maps involved
are the $\AW$, $\EML$ and $\SHI$ operators given in the
Eilenberg--Zilber Theorem \cite{EZ59}. This theorem states that
there is a chain homotopy equivalence $(\AW,\EML, \SHI)$ from the
normalized chain complex $C^{\scst N} (K\times L)$ of the
cartesian product  of $K$ and $L$ to the tensor product $C^{\scst
N}(K) \otimes C^{\scst N}(L)$ of the normalized chain complexes
$C^{\scst N}(K)$ and $C^{\scst N}(L)$. Whereas the number of
summands in the formula for $\AW$ grows linearly, the number of
summands in the formulae for $\EML$ and $\SHI$ grow exponentially,
then in order to define
  ``computable'' algebraic--combinatorial invariants, it seems that the right strategy
 is reduced to
  determine compositions of maps in which the morphism $\AW$ is involved. For example, the cup
  product on cohomology is essentially determined at the cochain level by the morphism $\AW$ and
  the diagonal map. All of this fits well with the results of Kristensen \cite{Kri63,Kri67}, where
a representation result for stable primary and secondary
cohomology operations in terms of cochain maps is given; and that
of Klaus \cite{Klaus1,Klaus2}, extending Kristensen's results to prove
that any cohomology operation module $p$ can be described in
terms of polynomials of coface operators at the cochain level.
 This approach is also corroborated in \cite{Rea96}, \cite{GR99} and \cite{casc02} where
 Steenrod squares,
 Steenrod  $k$th powers and Adem secondary cohomology operations
are seen at the cochain level
 essentially as compositions
 of the type
 \begin{equation}
H= \AW_{(p)} t_r \SHI_{(p)} t_{r-1} \cdots \SHI_{(p)} t_1 \SHI_{(p)}:
 C^{\scst N}(K^{\scst \ti p})\rightarrow C^{\scst N}(K)^{\scst \ot p}
\label{h}
\end{equation}
 where $t_i$ are  permutations of $p$ factors and $\AW_{(p)}$ and $\SHI_{(p)}$ are,
 respectively, the
 $\AW$ and $\SHI$ operators given by the Eilenberg--Zilber Theorem for $p$ simplicial sets.
 It is evident that
 an algorithm for computing
 these cohomology operations based on the previous formulation shows  extremely high computational
 costs.
  Because of this, a normalization of  compositions of face and degeneracy operators and
   a following step of the elimination of those summands of the normalized formula for $H$ with a factor
  having a degeneracy operator in its expression are done in order. This ``simplification''
  process allows to reach to a combinatorial description for $H$ having the minimum number of face
  operators involved.

 In this paper, we work with a general simplicial expression of type (\ref{h}), where the
 $t_i$ can be any  permutation.  We have developed a   software using  {\em Mathematica}
that deduces its  ``minimal'' simplicial  formulation. In particular, the solution to
this combinatorial problem provides  a way to design an efficient
 algorithm for computing any  Steenrod cohomology operation on  any cohomology class of
any degree. This work has been presented in \cite{aca02}.

The paper is organized as follows: In Section \ref{back} we
review the necessary theoretical background. In Section \ref{simplification}
we develop simplification techniques for obtaining an ``economical'' formulation
for operations of the type (\ref{h}).
Finally,
 Section \ref{example} is devoted to show an application
of
our method: an algorithm
for computing the Steenrod $k$th power
$P_{p}^k$ on the cohomology of any locally finite simplicial set is developed.

\section{Preliminaries}
\label{back}

In this section we introduce the notation and terminology used
throughout this paper. References for this material appear in
\cite{May67} and \cite{McL75}.

A {\em simplicial set} $K$ is a graded set indexed by the non--negative
integers together with {\em face} and {\em degeneracy operators}
$\partial_{i}:K_{q}\rightarrow K_{q-1}$ and $s_{i}:K_{q}\rightarrow K_{q+1}$, $0\leq
i\leq q$,  satisfying the following identities:
$$
\begin{array}{rll}
\mbox{\bf (i)} & \partial _{i}\partial _{j}=\partial _{j-1}\partial _{i}, &
i<j; \\
\mbox{\bf (ii)} & s_{i}s_{j}=s_{j+1}s_{i}, & i\leq j; \\
\mbox{\bf (iii)} & \partial _{i}s_{j}=s_{j-1}\partial _{i}, & i<j; \\
& \partial _{i}s_{j}=s_{j}\partial _{i-1}, & i>j+1; \\
& \partial _{j}s_{j}=1_{K_q}=\partial _{j+1}s_{j}. &
\end{array}$$
The elements of $K_{q}$ are called $q$--{\em simplices}. A simplex $x$ is
{\em degenerate} if $x=s_{i}(y)$ for some simplex $y$ and degeneracy operator $s_{i}$;
otherwise, $x$ is {\em non--degenerate}. Let $K$ and $L$ be two simplicial sets.
A map $f=\sum
f_{q}:\,K_q\rightarrow L_q$ of degree zero is a {\em simplicial map }if it
commutes with face and degeneracy operators, i.e.,
$f_{q}\partial _{i}=\partial _{i}f_{q+1}$ and  $f_{q}s_{i}=s_{i}f_{q-1}$.

The {\em cartesian product} $K\times L$ is a simplicial set whose simplices and
face and degeneracy operators are given by
$$\begin{array}{c}
(K\times L)_{q}=K_{q}\times L_{q}\,, \qquad
\partial _{i}(x,y)=(\partial _{i}x,\partial _{i}y)\,,\qquad
s_{i}(x,y)=(s_{i}x,s_{i}y).\end{array}
$$

Let $R$ be a commutative ring with identity $1\neq 0$.
The {\em chain complex} of a simplicial set $K$ with coefficients in $R$,
denoted by
 $C_{\ast }(K)$, is constructed as follows. Let $C_{n}(K)$
denote the free $R$--module on the set $K_{n}$. The face operators
$\partial_{i}$ linearly extend to module maps
$\partial_i:\,C_{n}(K)\rightarrow C_{n-1}(K)$. The alternating
sum $$\D d_n=\sum_{i=0}^n (-1)^{i}\partial _{i}: C_n(K)\ra
C_{n-1}(K)$$ is an $R$--module endomorphism of degree $-1$ such
that $d_nd_{n+1}$ is null for every $n\geq 0$; it is called the {\em differential} on
$C_{\ast}(K)$. The normalized chain complex $C^{\scst N}_*(K)$ is
defined by the quotient $$ C^{\scst N}_n(K) =
C_n(K)/s(C_{n-1}(K)),$$ where $s(C_{n-1}(K))$ denotes the free
$R$--module on the set of all the degenerate $n$--simplices of
$K$. Since we  always work with normalized chain complexes, we
simplify notation and write $C_*(K)$ instead of $C^{\scst
N}_*(K)$. $Z_{n}(K)=\ker d_{n}$ is the module of $n$--{\em cycles}
in $C_{\ast }(K)$; $B_{n}(K)=\mbox{Im }d_{n+1}$ is the module of
$n$--{\em boundaries} in $C_{\ast}(K);$ the quotient
$H_{n}(K)=Z_{n}(K)/B_{n}(K)$ is the $n$th {\em homology module}
of $K$. The homology class of a cycle $a\in Z_{n}(K)$ is denoted
by $[a]$.

Given  an abelian group $G$, form the abelian group
$$C^{n}(K;G)=Hom_{R}(C_{n}(K),G)$$ for each $n;$ the elements of
$C^{n}(K)$ are called the $n$--{\em cochains} of $C^{\ast
}(K;G)$. The differential $d$ on $C_{\ast }(K)$ induces a
codifferential $\delta: C^{\ast }(K;G)\ra C^{\ast +1}(K;G)$ of
degree $+1$
 via $\delta c=cd;$
the {\em cohomology} of $K$ is the family of abelian groups
$$
H^{n}(K;G)=\ker \delta^{n}/\mbox{Im }\delta^{n-1}.
$$
$B^{n}(K;G)=\mbox{Im }\delta^{n-1}$ is the module of $n$--{\em coboundaries}$;$
$Z^{n}(K;G)=\ker \delta^{n}$ is the module of  $n$--cocycles.
Furthermore, if $G$ is a ring, $H^{\ast }(K;G)$\ is an algebra with\ respect to the
{\em cup product}
$$
\smile :H^{i}(K;G)\otimes H^{j}(K;G)\rightarrow H^{i+j}(K;G)
$$
defined for $[c^{i}]\in H^{i}(K;G)$ and $[c^{j}]\in H^{j}(K;G)$
by $ [ c^{i}]\smile [ c^{j}]=[c^{i}\smile c^{j}]$,  where
\label{cupcero}
\begin{eqnarray*}
\left( c^{i}\smile c^{j}\right) (x)=\mu (c^{i}(\partial
_{i+1}\cdots
\partial _{i+j}x)\otimes c^{j}(\partial _{0}\cdots \partial _{i-1}x))
\end{eqnarray*}
for $x\in C_{i+j}(K);$ here $\mu$ is the multiplication on $G$.

Whenever two graded objects $x$ and
$y$ of degree $p$ and $q$ are interchanged
we apply the {\em Koszul's convention} and introduce the sign
$\left( -1\right) ^{pq}.$ The {\em tensor product} of chain
complexes $C_{\ast }(K)$ and $C_{\ast }(L)$ is the chain complex
$C_{\ast }(K)\otimes C_{\ast }(L)$ with differential $d_{\scst
C_{\ast }(K)\otimes C_{\ast }(L)}=d_{\scst C_{\ast }(K)}\otimes
1_{\scst C_{\ast }(L)}+1_{\scst C_{\ast }(K)}\otimes d_{\scst
C_{\ast }(L)}.$ Thus if $x_{p}\in C_{p}(K)$ and $y_{q}\in
C_{q}(L)$, an application of the Koszul convention gives
\begin{eqnarray*}
d_{\scst C_{\ast }(K)\ot C_{\ast }(L)}\left( x_{p}\ot
y_{q}\right) &=&( d_{\scst C_{\ast }(K)}\ot 1_{\scst C_{\ast
}(L)} +1_{\scst C_{\ast }(K)}\ot d_{\scst C_{\ast }(L)})
(x_{p}\ot y_{q}) \\
&=&d_{\scst C_{\ast }(K)}( x_{p}) \ot y_{q}+( -1)^{q} x_{p}\ot
d_{\scst C_{\ast }(L)}( y_{q})\,.
\end{eqnarray*}
A module homomorphism $f:C_{\ast }(K)\rightarrow C_{\ast }(L)$ of
degree zero such that $df=fd${\em \ }is a {\em chain map. } If
$f:C_{\ast }(K)\rightarrow C_{\ast }(L)$ and $g:C_{\ast
}(K^{\prime })\rightarrow C_{\ast }(L^{\prime })$ are chain maps,
so is $f\otimes g:C_{\ast }(K)\otimes C_{\ast }(K^{\prime
})\rightarrow C_{\ast }(L)\otimes C_{\ast }(L^{\prime })$.
Examples of chain maps are:
\begin{itemize}
\item The {\em diagonal map} $\Delta: C_*(K)\ra C_*(K^{\ti n})$
defined by $\Delta(x)=(x,\stackrel{\mbox{\scsz $n$ times}}\dots,x)$.
\item The {\em cyclic permutations}  $$t: C_*(K^{\ti n})\ra C_*(K^{\ti n})\qquad\mbox{ and }
\qquad T: C_*(K)^{\ot n}\ra C_*(K)^{\ot n}$$ such that
$$t(x_1,x_2,\dots,x_n)=(x_2,\dots,x_n,x_1)$$ and $$T(x_1\ot
x_2\ot \dots\ot x_n)=(-1)^{|x_1|(|x_2|+\cdots+|x_n|)}(x_2\ot
\dots\ot x_n\ot x_1)\,.$$
\end{itemize}
A {\em contraction from }$C_{\ast }(K)${\em \ to }$C_{\ast }(L) $
is a triple of homomorphisms $r=(f,g,\phi ),$ respectively
referred to as the {\em projection}, {\em inclusion} and {\em
homotopy operator}, with the following properties:
\begin{itemize}
\item $f:C_{\ast }(K)\rightarrow C_{\ast }(L)$ is a surjective chain map,
\item $g:C_{\ast }(L)\rightarrow C_{\ast }(K)$ is an injective chain map,
\item $\phi :C_{\ast }(K)\rightarrow C_{\ast +1}(K)$ is an endomorphism of degree
$+1$,
\item $d_{\scst C_*(K)}\phi +\phi d_{\scst C_*(K)}=1_{\scst C_{\ast }(K)}-gf$.
\end{itemize}
\noindent Furthermore, $f,$ $g$ and $\phi $ satisfy the following identities:
$$
 \phi g=0\,, \qquad  f\phi =0 \qquad\mbox{and}\qquad  \phi \phi =0.$$
A contraction will be denoted by $r=(f,g,\phi ):C_{\ast
}(K)\Rightarrow C_{\ast }(L)$. Two contractions $r=(f,\,g,\,\phi
):\,C_*(K)\Rightarrow C_*(L)$ and $r^{\prime }=(f^{\prime
},\,g^{\prime },\,\phi ^{\prime }):\,C_*(K^{\prime })\Rightarrow
C_*(L^{\prime })$ can be canonically combined to form new
contractions in the following ways:
\begin{itemize}
\item The {\em tensor product contraction} given by
$$
r\otimes r^{\prime }=(f\otimes f^{\prime },\,g\otimes g^{\prime },\,\phi
\otimes g^{\prime }f^{\prime }+1\otimes \phi ^{\prime }):C_*(K)\otimes
C_*(K^{\prime })\Rightarrow C_*(L)\otimes C_*(L^{\prime })\,.
$$
\item If $L=K'$, the {\em composition contraction} given by
$$
r^{\prime }r=(f^{\prime }f,\,gg^{\prime },\,\phi +g\phi ^{\prime
}f):C_*(K)\Rightarrow C_*(L^{\prime })\,.
$$
\end{itemize}

Let $p$ and $q$ be non--negative integers. A {\em $(p,q)$--shuffle} $(\alpha
,\beta )$ is a partition
$$\left\{ \alpha _{1}<\cdots <\alpha _{p}\right\}\;
\cup\; \left\{ \beta _{1}<\cdots <\beta _{q}\right\}$$ of the set $\left\{
0,1,\ldots ,p+q-1\right\} $. The signature of $(\alpha ,\beta )$ is given by
$$
sig(\alpha ,\beta )=\sum_{1\leq i\leq p}\alpha _{i}-(i-1).
$$
Let $\gamma =\{\gamma _{1},\dots ,\gamma _{r}\}$ be a set of
integers. Then
$s_{\gamma }$ denotes the
composition  of the degeneracy operators $s_{\gamma _{r}}\cdots
s_{\gamma _{1}}$.

An {\em Eilenberg--Zilber contraction} \cite{EZ59} from the chain
complex $C_{\ast }(K\times L)$ to the tensor product of chain
complexes $C_{\ast }(K)$ and $C_{\ast }(L)$ is a triple $r_{\scst
EZ}=(\AW,\EML,\SHI)$ such as:

\begin{itemize}
\item The Alexander--Whitney operator $\AW:C_{\ast}(K\times
L)\longrightarrow C_{\ast }(K)\otimes C_{\ast }(L)$ is defined by
$$
\AW(x_{m}, y_{m})=\sum_{0\leq i\leq m}\partial _{i+1}\cdots
\partial_{m}x_{m}\otimes \partial _{0}\cdots \partial _{i-1}y_{m}\,,
$$
where $(x_{m}, y_{m})\in C_m(K\times
L)$.
\item The Eilenberg--Mac Lane operator $\EML:C_{\ast }(K)\otimes C_{\ast
}(L)\longrightarrow C_{\ast }(K\times L)$ is defined by
$$
\EML(x_{p}\otimes y_{q})=\sum_{(\alpha ,\beta )\in \{(p,q)\mbox{\scsz --shuffles}
\}}(-1)^{sig(\alpha ,\beta )}(s_{\beta }x_{p},\,s_{\alpha }y_{q})\,,
$$
where $x_{p}\otimes y_{q}\in C_p(K)\otimes
C_q(L)$.
\item And the Shih operator $\SHI:C_{\ast }(K\times L)\longrightarrow
C_{\ast+1}(K\times L)$ is defined by
\begin{eqnarray*}
&&\SHI(x_{0},\,y_{0})=  0\,, \\
&&\SHI(x_{m},\,y_{m})\\
&&= \D \sum_{T(m)}\left( -1\right) ^{\epsilon(\al,\be) }(
s_{\bar{\beta}+\bar{m}}\partial _{m-q+1}\cdots \partial _{m}x_{m},s_{\alpha +
\bar{m}}\partial _{\bar{m}}\cdots \partial _{m-q-1}y_{m})\,,
\end{eqnarray*}
where
$$
\begin{array}{l}
T(m)=\{0\leq p\leq m-q-1\leq m-1\,,\;(\alpha ,\beta )\in \{(p+1,q)\mbox{--shuffles}\}\}, \\
\bar{m}=m-p-q, \\
\alpha +\bar{m}=\{\alpha _{1}+\bar{m},\dots ,\alpha _{p+1}+\bar{m}\}, \\
\bar{\beta}+\bar{m}=\{\bar{m}-1,\beta _{1}+\bar{m},\dots, \beta _{q}+\bar{m},\}, \\
\epsilon(\al,\be) =\bar{m}-1+sig(\alpha ,\beta ).
\end{array}
$$
\end{itemize}
 A recursive formula for the $\SHI$ operator appears in \cite{EM54}. The explicit formula
  given here was stated by Rubio in
 \cite{Rub91} and proved by Morace in the appendix of \cite{Rea96}.
It is evident that the $\AW$ operator has a
polynomial nature (concretely, the number
 of face  operators involved in its formula  is $O(m^2)$).
However, the $\EML$ and $\SHI$ operator have an essential
 ``exponential'' character, because
 shuffles  of degeneracy operators are involved in their respective formulations. In
  \cite{Prou83}, Prout\'{e}  determines that $\EML$ is unique and there is
  only two possibilities for $\AW$, both of its formulae being of the same complexity. Concerning
  $\SHI$, all the possible formulae have in common their exponential nature.

There is a contraction from $C_*(K^{\ti n})$ to
$C_*(K)^{\ot n}$ obtained by appropriately composing
Eilenberg--Zilber contractions. For any positive integers $s< n$,
let us denote by $r_{\scst EZ(n,s)}=(\AW_{(n,s)},\EML_{(n,s)},\SHI_{(n,s)})$
the contraction
$$r_{\scst
EZ(n,s)}\ot 1^{\ot s-1}=(\AW\ot 1^{\ot s-1},\EML\ot 1^{\ot
s-1},\SHI\ot 1^{\ot s-1})$$
from $C_*(K^{\ti n-s}\ti K)\ot C_*(K)^{\ot
s-1}$ to $C_*(K^{\ti n-s})\ot C_*(K)\ot C_*(K)^{\ot s-1}$.
Then, the composition $r_{\scst
EZ(n,n-1)}\cdots r_{\scst EZ(n,2)} r_{\scst EZ(n,1)}$ is a
contraction from $C_*(K^{\ti n})$ to $C_*(K)^{\ot n}$. We denote
it by $$r_{\scst EZ(n)}=(\AW_{(n)},\EML_{(n)},\SHI_{(n)}):C_*(K^{\ti
n})\Ra C_*(K)^{\ot n}\,.$$ Observe that the expression of
$\AW_{(n)}$ is:
\begin{eqnarray}\label{awn}
\begin{array}{cccl}
\AW_{(n)}({\bf x})&=&\AW_{(n,n-1)}&\AW_{(n,n-2)}\cdots\AW_{(n,2)}\AW_{(n,1)}
({\bf x})\\
&=& \D\sum_{0\leq i_1\cdots\leq i_{n-1}\leq m}&
\pa_{i_1+1}\cdots\pa_{m}x_1\\
&&&\ot \pa_0\cdots\pa_{i_1-1}\pa_{i_2+1}\cdots\pa_{m}x_2\\
&&&\vdots\\
&&&\ot \pa_0\cdots\pa_{i_{n-2}-1}\pa_{i_{n-1}+1}\cdots\pa_{m}x_{n-1}\\
&&&\ot \pa_0\cdots\pa_{i_{n-1}-1}x_{n}
\end{array}\end{eqnarray}
where ${\bf x}=(x_1,\dots x_n)\in C_m(K^{\ti n})$.
The number of face operators taking part in this formula is $O(n\cdot m^{n})$.

On the other hand, the expression of $\SHI_{(n)}$  in terms of the component morphisms of
the previous Eilenberg--Zilber contractions is:
\begin{eqnarray*}
\begin{array}{rl}
\D \sum_{\scst  1\leq \ell+1< n}&\EML_{\scst (n,1)}\cdots
\EML_{\scst (n,\ell)}\SHI_{\scst (n,\ell+1)}\AW_{\scst (n,\ell)}\cdots
\AW_{\scst (n,1)}\\
&=\SHI_{\scst (n,1)}\\
&\quad+\EML_{\scst (n,1)}\SHI_{\scst (n,2)}\AW_{\scst (n,1)}\\
&\quad\;\vdots\\
&\quad+\EML_{\scst (n,1)}\cdots \EML_{\scst (n,n-2)}\SHI_{\scst (n,n-1)}
\AW_{\scst (n,n-2)}\cdots \AW_{\scst (n,1)}.
\end{array}\end{eqnarray*}
Observe that whereas the number of summands in the formula for $\AW_{(n)}$ grows
in polynomial time (fixed $n$), the number of summands in the formulae for $\EML_{(n)}$
and $\SHI_{(n)}$ grow
exponentially.

\section{Simplification Techniques}\label{simplification}

 Let us recall that our motivation here is
to simplify  any composition of the type
$$\AW_{(p)}t_r\SHI_{(p)}t_{r-1}\cdots\SHI_{(p)} t_1 \SHI_{(p)}
=\sum \AW_{(p)}t_rESA_{(p,\ell_r)}\cdots t_1 ESA_{(p,\ell_1)}
$$
where every $t_i$ is any kind of permutation of $p$ factors,
$$ESA_{(p,\ell)}=\EML_{(p,1)}\cdots \EML_{(p,\ell)}
\SHI_{(p,\ell+1)}\AW_{(p,\ell)}\cdots \AW_{(p,1)}$$
and the sum is taken over the set
$\{1\leq i\leq r$, $ 0\leq \ell_i\leq p-2$, $1\leq k_i\leq p-1\}$.

We  will use the following basic properties:
\begin{itemize}
\item Any composition of face and
degeneracy operators of $K$ can be put in a
unique ``normalized" form: $$s_{j_t}\cdots
s_{j_1}\pa_{i_1}\cdots\pa_{i_r}\,,$$ where $j_t>\cdots>j_1\geq
0\;$ and  $\; i_s>\cdots >i_1\geq 0$.
\item Any summand on the tensor product of $n$ copies of $C_*(K)$ having a
factor (in the normalized form) with degeneracy operators in its
expression, is degenerate.
\end{itemize}

Let $i,j,m$ be integers such that $0\leq i\leq j\leq m$. The
interval $[i,j)$ denotes the set of consecutive integers from
$i$  to $j-1$.
\begin{itemize}
\item The {\em face--interval} $\pa_{[i,j)}$,
denotes the composition
$\pa_0\cdots\pa_{i-1}\pa_{j+1}\cdots\pa_m\,.$
\item If $i=0$ then $\partial_{[0,j)}=\pa_{j+1}\cdots\pa_m$.
\item If $j=m$ then
$\partial_{[i,m)}=\pa_0\cdots\pa_{i-1}$.
\item In the
case $i=j$ then
$\partial_{[i,i)}=\pa_0\cdots\pa_{i-1}\pa_{i+1}\cdots\pa_m$.
\end{itemize}
The notation $\partial_{[i,j)}$ must be interpreted as the
interval $[i,j)$ representing the indexes
 $\ell$, $0\leq\ell\leq m-1$,
such that $\pa_0\cdots\pa_{i_{j-1}-1} \pa_{i_j+1}\cdots\pa_m
s_{\ell}$ is degenerate. Whereas $j_1\leq i_2$, define the
following ``composition": $$\partial_{[i_1,j_1)}
\partial_{[i_2,j_2)}=\pa_0\cdots\pa_{i_1-1}\pa_{j_1+1}\cdots\pa_{i_2-1}\pa_{j_2+1}
\cdots\pa_m\,.$$ This composition can be extended without
problems to the composition of $n$ face--intervals.

With the new notation, we can rewrite the expression of
$\AW_{(n)}$ given in page \pageref{awn} as: $$\AW_{(n)}({\bf x})=\D
\sum_{P(m,n)} \partial_{[1]}x_1\ot
\partial_{[2]}x_2 \cdots\ot
\partial_{[n]}x_n\,,$$
where
 $[\ell]$ represents the interval $[i_{\ell-1},i_{\ell})$ and $P(m,n)$ is the set of all
 the possible partitions of $[0,m+1)$ in $n$ intervals.

First, in order to gradually show our technique, let us simplify
the composition $\AW_{(n)}t^k ESA_{(n,0)}({\bf x})=\AW_{(n)}t^k \SHI_{(n,1)}({\bf x})$, where $1\leq k\leq n-1$:
\begin{eqnarray}\label{0}\begin{array}{ll}
\AW_{(n)}t^k ESA_{(n,0)}({\bf x})&\\
=\D\sum_{ P(m+1,n)}\;\sum_{ T(m)}(-1)^{\epsilon(\al,\be)}&\partial_{[1]}s_{\be+\bar{m}}\pa_{m-q+1}\cdots\pa_mx_{k+1}\\
&\vdots\\
&\ot \partial_{[n-k-1]}s_{\be+\bar{m}}\pa_{m-q+1}\cdots\pa_mx_{n-1}\\
&\ot \partial_{[n-k]}s_{\al+\bar{m}}\pa_{\bar{m}}\cdots\pa_{m-q-1}x_n\\
&\ot \partial_{[n-k+1]}s_{\be+\bar{m}}\pa_{m-q+1}\cdots\pa_mx_1\\
&\vdots\\
&\ot
\partial_{[n]}s_{\be+\bar{m}}\pa_{m-q+1}\cdots\pa_mx_{k}\,.
\end{array}\end{eqnarray}
On one hand, $$(\alpha+\m)\cup(\beta+\m)=[\m-1,m+1)\quad\mbox{ and }\quad
\m-1\in\be+\m\,.$$ On the other hand, the non--degenerate summands
of (\ref{0}) satisfy that
$$(\alpha+\m)\cap [i_{n-k-1},i_{n-k})=\emptyset\quad\mbox{ and }\quad
(\beta+\m)\cap \left([0,i_{n-k-1})\cup
[i_{n-k},m+1)\right)=\emptyset\,.$$ We immediately obtain that
$$\be+\bar{m}\subset [i_{n-k-1},i_{n-k})\qquad \mbox{ and }\qquad
\al+\bar{m}\subset [0,i_{n-k-1})\cup [i_{n-k},m+1)\,,$$ therefore
$ i_{n-k-1}\leq \bar{m}-1$, $\;i_{n-k}=m-p$,
$$\be+\bar{m}=[\bar{m}-1,i_{n-k})\qquad\mbox{ and }\qquad
\al+\bar{m}=[i_{n-k},m+1)\,.$$ Now, we denote $$i'_j=\left\{
\begin{array}{cl}
i_j\quad&0\leq j<n-k\,,\\
i_j-q-1 \quad&n-k\leq j\leq n\,,\\
m\quad&j=n+1\,.
\end{array}
\right.$$

\newpage

\noindent and we can rewrite (\ref{0}) as:
\begin{eqnarray*}
\sum_{P(m,n+1)}(-1)^{ \tau_0}&&\partial_{ [1]}x_{ k+1} \ot \cdots
\ot\partial_{[n-k-1]}x_{n-1}
\ot \partial_{ [n-k]}\partial_{[n+1]}x_{n}\\
&&\ot\partial_{[n-k+1]}x_{1} \ot\cdots
\ot\partial_{[n]}x_{k}\end{eqnarray*} where
\begin{eqnarray*}
\tau_0&=
&\bar{m}-1+(p+1)q=i'_{n-k}+(i'_n-i'_{n-k})(i'_{n+1}-i'_n)\\&=&
|1|+\cdots+|n-k|+(|n-k+1|+\cdots +|n|)|n+1|\,,\end{eqnarray*} $|\ell|$ being
$i'_{\ell}-i'_{\ell-1}$.

In the same way, the expression of
$\AW_{(n)}t^{k} ESA_{(n,1)}({\bf x})$ is:
\begin{eqnarray}\label{1}\begin{array}{cl}
\D\sum_{\stackrel{\stackrel{\scst 0\leq \iota\leq m}{\scst P(m+1,n),\;  T(\iota)}}{
\scst (a,b)\in \{(\iota+1,m-\iota)\mbox{\scsz --sh.}\}}}&(-1)^{\scst sig(a,b)+\epsilon(\al,\be)}\\
&\partial_{[1]}s_b s_{\be+\bar{\iota}}\pa_{\iota-q+1}\cdots\pa_mx_{k+1}\\
&\vdots\\
&\ot \partial_{[n-k-2]}s_b s_{\be+\bar{\iota}}\pa_{\iota-q+1}\cdots\pa_mx_{n-2}\\
&\ot \partial_{[n-k-1]}s_b s_{\al+\bar{\iota}}\pa_{\bar{\iota}}\cdots\pa_{\iota-q-1}\pa_{\iota+1}\cdots\pa_mx_{n-1}\\
&\ot \partial_{[n-k]}s_{a}\pa_0\cdots\pa_{\iota-1}x_n\\
&\ot \partial_{[n-k+1]}s_b s_{\be+\bar{\iota}}\pa_{\iota-q+1}\cdots\pa_mx_1\\
&\vdots\\
&\ot
\partial_{[n]}s_b s_{\be+\bar{\iota}}\pa_{\iota-q+1}\cdots\pa_mx_{k}\,.
\end{array}\end{eqnarray}
On one hand, $a\cup b=[0,m+1)$ and on the other hand, the
non--degenerate summands satisfy that $$a\cap
[i_{n-k-1},i_{n-k})=\emptyset\qquad \mbox{ and }\qquad b\cap
([0,i_{n-k-1})\cup [i_{n-k},m+1))=\emptyset,$$ then
$b=[i_{n-k-1},i_{n-k})$ and $a=[0,i_{n-k-1})\cup [i_{n-k},m+1)$.
We denote $$i'_j=\left\{
\begin{array}{cl}
i_{j}\quad&0\leq j<n-k\,,\\
i_{j+1}-m+\iota\quad&n-k\leq j\leq n-2\,.
\end{array}
\right.$$

\newpage
\noindent Therefore (\ref{1}) becomes
\begin{eqnarray}\label{01}\begin{array}{ll}
\D\sum_{\stackrel{0\leq\iota\leq m}{\scst P(\iota+1,n-1),\;T(\iota)}}(-1)^{\scst sig(a,b)+\epsilon(\al,\be)}
&\partial_{[1]} s_{\be+\ov{\iota}}\pa_{\iota-q+1}\cdots\pa_{m}x_{k+1}\\
&\vdots\\
&\ot \partial_{[n-k-2]} s_{\be+\ov{\iota}}\pa_{\iota-q+1}\cdots\pa_{m}x_{n-2}\\
&\ot \partial_{[n-k-1]}
s_{\al+\ov{\iota}}\pa_{\bar{\iota}}\cdots\pa_{\iota-q-1}
\pa_{\iota+1}\cdots\pa_{m}x_{n-1}\\
&\ot \pa_{0}\cdots\pa_{\iota-1}x_n\\
&\ot \partial_{[n-k]} s_{\be+\ov{\iota}}\pa_{\iota-q+1}\cdots\pa_{m}x_1\\
&\vdots\\
&\ot
\partial_{[n-1]} s_{\be+\ov{\iota}}\pa_{\iota-q+1}\cdots\pa_{m}x_{k}\,.
\end{array}\end{eqnarray}
and  $sig(a,b)$ is $(m-\iota)(\iota+1-i'_{n-k-1})$. Now, we can
observe that if $k+1=n$ then the composition above is degenerate,
else $$i'_{n-k-2}\leq \bar{\iota}-1\,,\quad
\be+\bar{\iota}=[\bar{\iota}-1,i'_{n-k-2}) \quad\mbox{ and }\quad
\al+\bar{\iota}=[i'_{n-k-1},\iota+1)\,.$$ We denote
$$i''_j=\left\{
\begin{array}{cl}
i'_j\quad& 0\leq j<n-k-1\,,\\
i'_{ j}-q-1\quad & n-k-1\leq j\leq n-2\,,\\
\iota-q\quad & j=n-1\,,\\
\iota\quad & j=n\,,\\
m\quad & j=n+1
\end{array}
\right.$$
therefore (\ref{01}) is
\begin{eqnarray*}
\D\sum_{P(m,n+1)}(-1)^{\tau_1} &&\partial_{[1]}x_{k+1} \ot\cdots
\ot \partial_{[n-k-2]} x_{n-2} \ot
\partial_{[n-k-1]}\pa_{[n]}x_{n-1}
\ot \pa_{[n+1]}x_n\\
&&\ot \partial_{[n-k]}x_1 \ot\cdots \ot\partial_{[n-1]}
x_{k}\end{eqnarray*} and the sign:
\begin{eqnarray*}
\tau_1&=&(m-\iota)(\iota+1-i'_{n-k-1})+\bar{\iota}-1+(p+1)q\\
&=&(i''_{n+1}-i''_n)(i''_n+1-i''_{n-k-1}+i''_{n-1}-i''_{n}-1)\\
&&+i''_{n-k-1}+(i''_{n-1}-i''_{n-k-1})(i''_n-i''_{n-1})\\
&=&i''_{n-k-1}+(i''_{n-1}-i''_{n-k-1})(i''_{n+1}-i''_{n-1})\\
&=&|1|+\cdots+|n-k-1|+(|n-k|+\cdots +|n-1|)(|n|+|n+1|)\,.\end{eqnarray*}

Now, let us study the general case.  As we said before, we are interested in simplifying any composition of the form
\begin{eqnarray}
\label{goal}
\AW_{(n)}t_r ESA_{(n,\ell_r)}\cdots t_1ESA_{(n,\ell_1)}\,.
\end{eqnarray}
We will do it inductively.
Let $h: C_{*}(K^{\ti n})\ra C_{*}(K)^{\ot n}$ be a
morphism of degree $r$ whose normalized expression is:
\begin{eqnarray*}
h({\bf x})=\sum_{P(m,n+r)} (-1)^{\scst sign\{[1],\dots,[n+r]\}}
\pa_{[\;]}x_{k_1}\ot\pa_{[\;]}x_{k_2}\cdots  \ot
\pa_{[\;]}x_{k_n}
\end{eqnarray*}
such that
$(x_{k_1},\dots,x_{k_n})=t_{\lambda}(x_1,\dots,x_n)$ where $t_{\lambda}:C_*(K^{\ti n})\ra C_*(K^{\ti n})$
is any permutation.
and each $\pa_{[\;]}$ denotes a composition of
non--consecutive elements of the set
$\{\pa_{[1]},\pa_{[2]},\dots,\pa_{[n+r]}\}$ where
$\{[1],[2],\dots,[n+r]\}\in P(m,n+r)$; moreover, each
$\pa_{[j]}$, $1\leq j\leq n+r$, appears exactly once in the
expression of $h({\bf x})$. Our goal is to simplify  the
composition $H=h\, ESA_{(n,\ell)}$,
 where $0\leq \ell\leq n-2$.
\begin{prop}\label{nodegenerado} If one of the following conditions holds on $h$:
\begin{itemize}
\item There is not any face--interval preceding $x_j$ for $1\leq j\leq n$;
\item There exists a factor in $h({\bf x})$ with more than one face--interval preceding $x_{n+1-u}$
for some $1\leq u\leq \ell$;
\item The
face--interval $\pa_{[j]}$ immediately before $x_{n-\ell}$  in
$h({\bf x})$
satisfies that\\
$j=\mbox{ max }\{v\mbox{ such that $\pa_{[v]}$ appears preceding
some $x_u$ for $1\leq u\leq n-\ell$}\}\,;$
\end{itemize}
then all the summands of $H$ are degenerate.
\end{prop}

From now on, let us suppose that $h({\bf x})$ does not satisfy any of the conditions of the
proposition above.
Let us denote by $\pa_{[j_u]}$ the unique face--interval preceding $x_{n+1-u}$
for $1\leq u\leq \ell$.

\begin{lem}\label{nodegenerado2} If
  the composition $\pa_{[j_u-1]}\pa_{[j_u+1]}$ appears in the
expression of $h$
for some $u$, $1\leq u\leq \ell$,
 then all the summands of $H$ are degenerate.
\end{lem}

\begin{thm}\label{alg}{\sc Simplification Algorithm.}

\begin{tabbing}
{\sc Input:} \= {\tt The morphism $h:\,C_{*}(K^{\ti n})\ra
C_{*}(K)^{\ot n}$ of degree $r$ described}\\
\> {\tt  above such that it does
not  satisfy either  Proposition \ref{nodegenerado}}\\
\> {\tt  or Lemma
\ref{nodegenerado2}.}\\
 {\sc Output:} {\tt  The simplified expression of
$\;H({\bf x)}=h\,ESA_{(n,\ell)}({\bf x)}$.}
\end{tabbing}
\begin{tabbing}
{\tt For }\=  $u=1$ {\tt to }$u=\ell$ {\tt do}\\
\>{\tt replace  $\pa_{[j_u]}$ preceding $x_{n+1-u}$ by
$\pa_{[n+r+2-u]}$.}\\
{\tt End for.}
\end{tabbing}
Let $\{\pa_{[v_1]},\dots,\pa_{[v_{n+r-\ell}]}\}$,
$v_1<\cdots<v_{n+r-\ell}$, denote the set of the face-intervals
preceding $x_u$ for $1\leq u\leq n-\ell$.
\begin{tabbing}
{\tt For }\=  $s=1$ {\tt to }$s=n+r-\ell$ {\tt do}\\
\>{\tt replace  $\pa_{[v_s]}$  by
$\pa_{[s]}$.}\\
{\tt End for.}\\
{\tt Replace $x_{n-\ell}$ by $\pa_{[n+r-\ell+1]}x_{n-\ell}$.}
\end{tabbing}
\noindent Starting from the sign of $h$ of degree $m+1$, we obtain the sign of $H$ of
 degree $m$
 as follows.

 \begin{itemize}
 \item[] Step 1:
 \begin{tabbing}
{\tt  For }\= $u=1$ {\tt to }$u=\ell$ {\tt do}\\
\>{\tt replace $|j_u|$  by $|n+r-u+1|+1$.}\\
\>{\tt For }\= {\tt $j=j_u+1$ to $j=n+r-u+1$ do}\\
\>\>{\tt replace $|j|$ by $|j-1|$.}\\
\>{\tt End for;}\\
\>{\tt add $(|n+r-u+1|+1)(|j_u|+\cdots +|n+r-u|)$}\\
{\tt End for.}
\end{tabbing}
\item[]
\item[] Let $\pa_{[v]}$ be the face--interval immediately
before $x_{n-\ell}$. Starting from the modified sign of $H$ do
\item[]
\item[] Step 2:
\begin{tabbing}
{\tt For }\= {\tt $j=n+r-\ell+2$ to $j=n+r$ do}\\
\> {\tt replace $|j|$ by $|j+1|$.}\\
{\tt End for;}\\
{\tt replace $|n+r-\ell+1|$ by $|n+r-\ell+2|-1$;}\\
{\tt replace $|v|$ by $|n+r-\ell+1|+1$;}\\
{\tt add $|1|+\cdots +|v|+ (|v+1|+\cdots+|n+r-\ell|)|n+r-\ell+1|$.}\\
\end{tabbing}
\end{itemize}
\end{thm}

\Elproofname

 For the sake of simplicity but without lost of generality, we
consider that the expression of $h({\bf x)}$ is $$\sum_{P(m,n+r)}
(-1)^{\scst
sign\{[1],\dots,[n+r]\}}\pa_{[\;]}x_1\ot\cdots\pa_{[\;]}x_{n-\ell}\ot
\ot\pa_{[j_{\ell}]}x_{n-\ell+1}\ot\cdots\ot \pa_{[j_{1}]}x_n\,;$$
 consequently, the expression of $H({\bf x})$ is:
\begin{eqnarray}\label{24}\begin{array}{cl}
\D\sum_{
\stackrel{\scst P(m+1,n+r),\; T(\iota_{\ell})}{
\stackrel{\scst 0\leq \iota_{\ell}\leq \iota_{\ell-1}\leq \cdots \leq \iota_{1}\leq m}{
\scst \{(a_j,b_j)\in\{(\iota_j+1,m-\iota_j)\mbox{\scsz --sh.}\}:\; 1\leq j\leq \ell\}
}
}
}
&(-1)^{\scst sign\{[1],\dots,[n+r]\}+ sig(a_1,b_1)+\cdots+sig(a_{\ell},b_{\ell})+\epsilon(\al,\be)}\\
&\partial_{[\;]}s_{b_1}\cdots s_{b_{\ell}} s_{\be+\ov{\iota_{\ell}}}\pa_{\iota_{\ell}-q+1}\cdots\pa_mx_{1}\\
&\vdots\\
&\ot  \partial_{[\;]}s_{b_1}\cdots s_{b_{\ell}} s_{\be+\ov{\iota_{\ell}}}\pa_{\iota_{\ell}-q+1}\cdots\pa_mx_{n-\ell-1}\\
&\ot \partial_{[\;]}s_{b_1}\cdots s_{b_{\ell}} s_{\al+\ov{\iota_{\ell}}}\pa_{\ov{\iota_{\ell}}}\cdots\pa_{\iota_{\ell}-q-1}\pa_{\iota_{\ell}+1}\cdots\pa_mx_{n-\ell}\\
&\ot \partial_{[j_{\ell}]}s_{b_1}\cdots s_{b_{\ell-1}}s_{a_{\ell}} \pa_0\cdots\pa_{\iota_{\ell}-1}\pa_{\iota_{\ell-1}+1}\cdots\pa_mx_{n-\ell+1}\\
&\vdots\\
&\ot \partial_{[j_2]}s_{b_1}s_{a_2} \pa_0\cdots\pa_{\iota_2-1}\pa_{\iota_1+1}\cdots\pa_mx_{n-1}\\
&\ot \partial_{[j_1]}s_{a_1}\pa_0\cdots\pa_{\iota_1-1}x_n\,.
\end{array}\end{eqnarray}
The non--degenerate summands of $H({\bf x})$ satisfy that
$$a_1=[0,i_{j_1-1})\cup [i_{j_1},m+1)\qquad\mbox{ and }\qquad
b_1=[i_{j_1-1},i_{j_1})\,.$$ Then, $$\begin{array}{l}
i^1_j=i_j\,\mbox{ for }0\leq j<j_1,\\
i^1_j=i_{j+1}-m+\iota_1\,\mbox{ for }j_1\leq j< n+r-1,\\
i^1_{n+r-1}=\iota_1+1,\\i^1_{n+r}=m.\end{array}$$
Therefore,  we have that
$$\begin{array}{l}
i_j=i^1_j\,\mbox{ for }0\leq j<j_1,\\
i_{j}=i^1_{j-1}+i^1_{n+r}-i^1_{n+r-1}+1
\,\mbox{ for }j_1\leq j\leq n+r-1,\\
i_{n+r}=i^1_{n+r}+1\,.\end{array}$$ So, in
$sign\{[1],\dots,[n+r]\}$, $|j_1|$ is replaced by $|n+r|+1$,
$\;|j|$ is replaced by $|j-1|$ for $j_1< j\leq n+r$ and
\begin{eqnarray*}
sig(a_1,b_1)&=&(m-\iota_1)(\iota_1+1-i_{j_1-1})
=(i^1_{n+r}-i^1_{n+r-1}+1)(i^1_{n+r-1}-i'_{j_1-1})\\
&=&(|n+r|+1)(|j_1|+\cdots+|n+r-1|)\,,
\end{eqnarray*}
is added.

In general, fixed $u$, $1\leq u\leq \ell$, we have that
$$a_u=[0,i^{u-1}_{j_u-1})\cup [i^{u-1}_{j_u},\iota_{u-1}+1)
\qquad\mbox{ and }\qquad b_u=[i^{u-1}_{j_u-1},i^{u-1}_{j_u})\,.$$
Then, $$\begin{array}{l}
i^u_j=i^{u-1}_j\,\mbox{ for }0\leq j<j_u,\\
i^u_j=i^{u-1}_{j+1}-\iota_{u-1}+\iota_u\,\mbox{ for }j_u\leq j< n+r-u,\\
i^u_{n+r-u}=\iota_u+1,\\i^u_{n+r-u+1}=\iota_{u-1}.\end{array}$$
Therefore,
$$\begin{array}{l}
i^{u-1}_j=i^u_j\;\mbox{ for }\,0\leq j<j_u\,\mbox{ and }\,n+r-u+2\leq j\leq n+r,\\
i^{u-1}_{j}=i^u_{j-1}+i^u_{n+r-u+1}-i^u_{n+r-u}+1
\,\mbox{ for }j_u\leq j\leq n+r-u,\\
i^{u-1}_{n+r-u+1}=i^u_{n+r-u+1}+1\,.\end{array}$$ So, in
$sign\{[1],\dots,[n+r]\}$, $|j_u|$ is replaced by $|n+r-u+1|+1$
and $|j|$ is replaced by $|j-1|$ for $j_u< j\leq n+r-u+1$. Also,
\begin{eqnarray*}
sig(a_u,b_u)&=&(\iota_{u-1}-\iota_u)(\iota_u+1-i^{u-1}_{j_1-1})\\
&=&(|n+r-u+1|+1)(|j_u|+\cdots+|n+r-u|)\,,\end{eqnarray*}
is added.
Therefore, the expression of (\ref{24}) is:
\begin{eqnarray*}\begin{array}{cl}
\D\sum_{
\stackrel{\scst P(\iota_{\ell}+1,n+r-\ell),\; T(\iota_{\ell})}{
\scst 0\leq \iota_{\ell}\leq \iota_{\ell-1}\leq \cdots\leq \iota_{1}\leq m}
}
&(-1)^{sign\{[1],\dots,[n+r]\}+\epsilon(\al,\be)}\\
&\partial_{[\;]}s_{\be+\ov{\iota_{\ell}}}\pa_{\iota_{\ell}-q+1}\cdots\pa_mx_{1}\\
&\vdots\\
&\ot  \partial_{[\;]}s_{\be+\ov{\iota_{\ell}}}\pa_{\iota_{\ell}-q+1}\cdots\pa_mx_{n-\ell-1}\\
&\ot \cdots\partial_{[v]}s_{\al+\ov{\iota_{\ell}}}\pa_{\ov{\iota_{\ell}}}\cdots\pa_{\iota_{\ell}-q-1}\pa_{\iota_{\ell}+1}\cdots\pa_mx_{n-\ell}\\
&\ot  \pa_0\cdots\pa_{\iota_{\ell}-1}\pa_{\iota_{\ell-1}+1}\cdots\pa_mx_{n-\ell+1}\\
&\vdots\\
&\ot \pa_0\cdots\pa_{\iota_2-1}\pa_{\iota_1+1}\cdots\pa_mx_{n-1}\\
&\ot \pa_0\cdots\pa_{\iota_1-1}x_n\,.
\end{array}\end{eqnarray*}

Now, $\al+\ov{\iota_{\ell}}=[i^{\ell}_{v},\iota_{\ell}+1)$ and
$\be+\ov{\iota_{\ell}}=[\ov{\iota_{\ell}}-1,i^{\ell}_{v})$, then
$$ i^{\ell+1}_j=i^{\ell}_j \;\mbox{ for }0\leq j\leq v-1\,,\qquad
i^{\ell+1}_j=i^{\ell}_j-q-1\;\mbox{ for }v\leq j\leq
n+r-\ell-1\,,$$ $$i^{\ell+1}_{n+r-\ell}=\iota_{\ell}-q,\qquad
i^{\ell+1}_{n+r-\ell+1}=\iota_{\ell}\,,\qquad
i^{\ell+1}_{j+1}=i^{\ell}_{j}\;\mbox{ for }n+r-\ell+1\leq j\leq
n+r\,.$$ That is, $$ i^{\ell}_j=i_j^{\ell+1} \;\mbox{ for }0\leq
j\leq v-1\,,\qquad i^{\ell}_j=i^{\ell+1}_j+q+1\,\mbox{ for }v\leq
j\leq n+r-\ell-1\,,$$
$$i^{\ell}_{n+r-\ell}=i^{\ell+1}_{n+r-\ell+1}+1\,,\qquad
i^{\ell}_{j}=i^{\ell+1}_{j+1}\,\mbox{ for }n+r-\ell+1\leq j\leq
n+r\,.$$ So, in $sign\{[1],\dots,[n+r]\}$, $|j|$ is replaced by
$|j+1|$ for $n+r-\ell+2\leq j\leq n+r$, $|v|$ is replaced by
$|n+r-\ell+1|+1$ and
 $|n+r-\ell+1|$ is replaced by $|n+r-\ell+2|-1$. Finally,
 \begin{eqnarray*}
\epsilon(\al,\be)&=&\ov{\iota_{\ell}}-1+(p+1)q\\
&=&i_v^{\ell+1}+(i_{n+r-\ell}^{\ell+1}-i_v^{\ell+1})
(i_{n+r-\ell+1}^{\ell+1}-i_{n+r-\ell}^{\ell+1})\\
 &=&|1|+\cdots+|v|+(|v+1|+\cdots+|n+r-\ell|)|n+r-\ell+1|
 \end{eqnarray*}
 is added.
\qed

\begin{thm}
The number of face operators taking part in the normalized formula for
$\AW_{(p)}t_r\SHI_{(p)}\cdots t_1 \SHI_{(p)}$ is,
in the worst case, $O(p^{r+1} m^{p+r+1})$.
\end{thm}

\Elproofname

On one hand, the number of summands of the form (\ref{goal}) is
$(p-1)^r$.
On the other hand, the number of summands in the simplified
formula for each morphism (\ref{goal}) is $O(m^{p+r})$ and the number of face
operators in each summand is $O(pm)$ .
Therefore the number of face operators taking part in the
normalized formula for $\AW_{(p)}t_r\SHI_{(p)}t_{r-1}\cdots\SHI_{(p)} t_1 \SHI_{(p)}$
 is
 $O((p-1)^rm^{p+r}pm)$ that is $O(p^{r+1}m^{p+r+1})$. \qed

\section{An Example: Algorithm for Computing $P_p^k$}\label{example}

In this section we study  the computation of the   cohomology
operations  Steenrod $k$th powers $P_p^k$ \cite{Ste52} as an
application of the technique given in the section above. First,
we give the definition of these operations at the cochain level
due to Steenrod \cite{Ste52}. We next  show explicit formulae developed in
\cite{GR99} for these operations in terms of Eilenberg--Zilber
contractions at the cochain level. Finally, we develop an
algorithm for computing $P_p^k$ at the cohomology level on any
locally finite simplicial set.

An infinite sequence of morphisms
$\{D^n_r:\,C_*(K)\ra C_*(K)^{\ot n}\}_{r\geq 0}$ of degree $r$ such that:
\begin{eqnarray}
\label{rel}D^n_0=\AW_{(n)}\Delta\,; \qquad d_{\scst C_*(K)^{\ot n}}
D^n_r+(-1)^{r-1}D^n_rd_{\scst C_*(K)}=\alpha_rD^n_{r-1}\,,\;r>0;\end{eqnarray}
where $\alpha_r:\,C_*(K)^{\ot n}\ra C_*(K)^{\ot n}$ is defined by
$$\alpha_r=\left\{\begin{array}{cl}
T-1&\mbox{ if $r$ odd,}\\
1+T+\cdots+T^{n-1}\quad&\mbox{ if $r$ even,}
\end{array}\right.$$
called a {\em higher diagonal approximation} \cite{Ste52}
``measures" the lack of commutativity of
 $\AW_{(n)}$.

In the particular case of
$p= 2$, it is possible to define cochain mappings called
 {\em cup--$i$ product}, $$\smile_i:C^q(K;G)\ot C^p(K;G)\ra C^{q+p-i}(K;G)$$  by
$c\smile_i c'=\mu (c\ot c')D^2_i\,.$ Observe that the expression
of $c\smile_0 c'$ coincides with that of the cup product given in
page \pageref{cupcero}. Taking $[c]\in H^j(K;\Zr)$, the
cohomology operations {\em Steenrod squares} \cite{Ste47} are
defined by
  $Sq^{i}[c]=[c\smile_{j-i} c]\in H^{j+i}(K;\Zr)$.

Now, let  $p>2$ be a prime number. Starting from the sequence (\ref{rel}), the
{\em Steenrod  $k$th power}
$P_p^k: H^q(K;{\bf Z}_p)\ra H^{q+2k(p-1)}(K;{\bf Z}_p)$, $q\geq 2k$,
 is defined at the cochain level as follows. If $c\in Z^q(K;{\bf Z}_p)$, then
\begin{eqnarray}\label{popo}
P^k_p(c)=R\,\mu c^{\ot p}
D^p_{(q-2k)(p-1)}\in Z^{q+2k(p-1)}(K;{\bf Z}_p)\,,\end{eqnarray} where
$\mu$ is the natural product on ${\bf Z}_p$ and
$R=(-1)^{(p-1)(k+\frac{1}{2}q(q-1))}\left(\left(\frac{p-1}{2}\right)!\right)^{2k-q}$.

The acyclic model method \cite{EM53} is used for guaranteeing the
existence of the morphisms $D^n_r$ ($n$ and $r$ being
non--negative integers). An alternative of the previous method is
to obtain the morphisms $D^n_r$ using algebraic fibrations with a
cartesian product of $n$ copies of a given simplicial set $K$ as
the base space and a subgroup of the symmetric group $S_n$ as the
fiber space. This last point of view has been established in
\cite{Rea96} and  \cite{GR99} for Steenrod operations, in
\cite{casc02} for secondary cohomology operations and generalized
in \cite{Gon00} for any cohomology operation. In \cite{GR99} we
obtain  explicit formulae for a higher diagonal approximation in
terms of the component morphisms of a given Eilenberg--Zilber
contraction. Let $\gamma_j:\,C_*(K^{\ti n})\ra C_*(K^{\ti n})$
define by $$\gamma_j=\left\{\begin{array}{cl}
t\quad&\mbox{ if $j$ odd}\\
t+\cdots+t^{n-1}\quad&\mbox{ if $j$ even.}
\end{array}\right.$$
then
$$D^n_r=\AW_{(n)}\gamma_r\SHI_{(n)}\cdots\gamma_1\SHI_{(n)}\Delta
=\D \sum\;  \AW_{(n)}t^{k_r}ESA_{(n,\ell_r)}\cdots t^{k_1}ESA_{(n\ell_1)}\Delta$$
where the sum is taken over all the possible
$1\leq \ell_i+1,k_i< n$,
where $k_i=1$ if $i+r$ odd;
for all $1\leq i\leq r$.

Observe that an algorithm based on these formulae for $D_r^n$ is not useful in practice,
due to the exponential nature of the morphisms involved. Nevertheless, we can
apply the Simplification Algorithm explained before in order to obtain a pure combinatorial definition
of $D_r^n$ only in terms of face operators.
Notice that
for obtaining a normalized expression of $D_r^n$, we have to apply Theorem \ref{alg}
$(n-1)^{\lfloor r/2\rfloor}(n-1)^r$ times in the worst case. However,
taking into account Proposition \ref{nodegenerado},
the non--degenerate summands of $D_r^n$ can only appear
when  $k_i+\ell_i<n$ for $1\leq i\leq r$. Moreover,  if $k_i+\ell_i<n$ and $k_i< \ell_{i+1}$
then the non--degenerate summands of $D_r^n$ can only appear
when $k_i+\ell_i< \ell_{i+1}$ for $1\leq i< r$.
Examples of the simplification process are:
\begin{eqnarray*}
D_1^n(x)=\sum_{P(m,n+1)}(-1)^{\tau_1}&&
\pa_{[1]}x\ot\cdots\ot \pa_{[n-\ell-2]}x\ot
\pa_{[n-\ell-1]}\pa_{[n-\ell+1]} x\\&&\ot\pa_{[n-\ell+2]}x\ot
\cdots\ot\pa_{[n+1]}x \ot\pa_{[n-\ell]}x\,,\end{eqnarray*}
where $\tau_1=|1|+\cdots+|n-\ell-1|+|n-\ell|(|n-\ell+1|+\cdots+|n+1|)$
and
\begin{eqnarray*}
D_2^n(x)
=
\sum_{\stackrel{\stackrel{\scst 0<\ell_2+1 \leq \ell_1<n-1}{\scst 0<k<n}}
{\scst P(m,n+2)}}(-1)^{\tau_2
}
&&\pa_{[1]}x\ot\cdots\ot\pa_{[n-k-\ell_1-1]}x
\ot\pa_{[n-k-\ell_1]}\pa_{[n-\ell_1+1]}x\\
&&\ot \pa_{[n-\ell_1+2]}x\ot\cdots\ot
\pa_{[n-\ell_2-1]}x\ot\pa_{[n-\ell_2]}\pa_{[n-\ell_2+2]}x\\
&&\ot\pa_{[n-\ell_2+3]}x\ot\cdots\ot
\pa_{[n+2]}x\ot\pa_{[n-\ell_2+1]}x\\&&\ot \pa_{[n-k_1-\ell_1+1]}x\ot\cdots\ot\pa_{[n-\ell_1]}
x\\\\
- \sum_{\stackrel{\scst 0<\ell+1,k<n}{\scst P(m,n+2)}}(-1)^{\tau_3}
 &&\pa_{[1]}x\ot\cdots\ot\pa_{[n-k-\ell-2]}x
\ot\pa_{[n-k-\ell-1]}\pa_{[n-\ell+2]}x\\
&&\ot
\pa_{[n-\ell+3]}x\ot\cdots\ot\pa_{[n+2]}x_n \ot
\pa_{[n-k-\ell]}\pa_{[n-\ell+1]}x\\&&
\ot\pa_{[n-k-\ell+1]}x\ot\cdots\ot \pa_{[n-\ell]}x
\end{eqnarray*}
where $\tau_2=(|n-k_1-\ell_1+1|+\cdots+|n-\ell_1|)(|n-\ell_1+1|+\cdots+|n-\ell_2-1|$
$+|n-\ell_2+1|+1+\cdots
+|n+1|+1)+|n-k_1-\ell_1+1|+\cdots+|n-\ell_2|$\\
$+|n-\ell_2+1|(|n-\ell_2+2|+1+\cdots
+|n+2|+1)$\\
and $\tau_3=|n-k-\ell-1|+(|n-k-\ell+1|+|n-\ell+1|)(|n-\ell+2|+1)+|n-\ell+3|+1+\cdots+|n+2|+1+
(|n-k-\ell+1|+\cdots+|n-\ell|)(|n-\ell+1|+|n-\ell+4|+\cdots+|n+1|)$.

Taking into account the sign and organization of the intervals in a general summand
of the normalized  expression of $D_1^n$ and $D_2^n$, it should be possible to obtain a general
expression of any $D^n_r$ but this study exceeds the scope of this paper.

\newpage

On the other hand, bearing in mind the expression  at the cochain
level of  the Steenrod  power operation $P^k_p(c)$
 where $c\in Z^q(K,{\bf Z}_p)$,  since
 $c$ is a $q$--cochain, we only consider
those summands in the normalized formula for $D_{(q-2k)(p-1)}^p$
with exactly $2k(p-1)$ face
operators in each factor.

Since the explicit formulae for the Steenrod powers operations
$P_p^k$ are given at the cochain level, in order to design an
algorithm for computing them at the cohomology level, we first
compute an explicit contraction $(f,g,\phi)$ from $C_*(K)$ to
$H_*(K)$, $K$ being a simplicial set finite in each degree and
${\bf Z}_p$ being the ground ring. This contraction can be
constructed using the classical matrix algorithm
\cite{Mun84} based on reducing certain matrices (corresponding to
the differential at each degree) to their Smith normal form
\cite{GR01}. The complexity of this method is
$O(M^3)$ where $M$ is the number of simplices of $K$.

Since the ground ring is a field, then the homology and cohomology are isomorphic.
Moreover, if $\alpha$ is a generator of homology of degree $q$, then $\alpha^*: H_q(K)\ra {\bf Z_p}$
such that
$$\alpha^*(\beta)=\left\{\begin{array}{l}
0\quad\mbox{ if $ \alpha\neq\beta\in H_q(K)$}\\
1\quad\mbox{ if $\beta= \alpha$,}\\
\end{array}\right.$$
is a generator of cohomology of degree $q$. Fixed $k$, suppose that the
normalized description of the morphism $D^p_{(q-2k)(p-1)}$
obtained using Theorem \ref{alg}, and a contraction $(f,g,\phi)$
from $C_*(K)$ to $H_*(K)$  using
the algorithm described above are given.
 Then,
(\ref{popo}) becomes at the cohomology level as:
$$P_p^k(\alpha^*)=\D \sum_{j=1}^{u}
R \left(\mu (\alpha^*f)^{\ot p}
D^p_{(q-2k)(p-1)}g(\gamma_j)\right)\cdot \gamma_j^*$$
where $\{\gamma_1,\dots,\gamma_u\}$ is a basis of $H_{q+2k(p-1)}$.

Summing up, we have designed an algorithm for computing any Steenrod reduced $k$th powers
on any  class  of cohomology for any locally finite simplicial set.

\end{document}